\documentclass[12pt]{article}
\usepackage{amsfonts}
\usepackage{mathrsfs}
\usepackage{amsmath,amssymb}
\numberwithin{equation}{section}

\def\BB{\mathcal{B}}

\def\NN{\mathcal{N}}
\def\HH{\mathcal{H}}
\def\l{\lambda}
\def\LL{\mathcal{L}^s_{\lambda,\mu}}
\def\gg{\mathfrak{g}}

\def\vp{\varphi}

\def\sm{\!\setminus\!}
\def\Inv{\mathrm{Inv}(\mathcal {L}^s_{\l,\mu},\C)}

\def\cl{\centerline}
\def\ni{\noindent}
\def\rar{\longrightarrow}
\def\vs{\vspace*}
\def\C{\mathbb{C}}

\def\Z{\mathbb{Z}}

\def\DU{\!\diagup\!}

\def\QED{\hfill$\Box$\par}

\newtheorem{theo}{Theorem}[section]
\newtheorem{coro}[theo]{Corollary}

\newtheorem{lemm}[theo]{Lemma}
\newtheorem{prop}[theo]{Proposition}

\begin{document}

\baselineskip 18pt

\cl{{\bf Non-degenerate Symmetric Invariant Bilinear Forms}}

\cl{{\bf on the Deformative Schr\"{o}dinger-Virasoro Algebras}
\footnote{Supported by NSF grants 10825101, 11101056, 11271056 of China}}

\vs{6pt}

\cl{Huanxia Fa$^{1,2)}$,\ \ Junbo Li$^{2)}$,\ \ Linsheng Zhu$^{2)}$}

\cl{\small $^{1)}$Wu Wen-Tsun Key Laboratory of Mathematics,}\vs{-3pt}
\cl{\small University of Science and Technology of China, Hefei 230026, China}

\cl{\small $^{2)}$School of Mathematics and Statistics, Changshu Institute
of Technology, Changshu 215500, China}

\cl{\small E-mail: sd\_huanxia@163.com,\ \ sd\_junbo@163.com,\ \ lszhu@cslg.edu.cn}\vs{6pt}

\noindent{\small{\bf Abstract.} In the present paper
we shall determine all the non-degenerate symmetric invariant bilinear forms on the deformative Schr\"{o}dinger-Virasoro algebras.

}


\vs{12pt}

\cl{\bf\S1. \
Introduction}\setcounter{section}{1}\setcounter{theo}{0}\setcounter{equation}{0}

The Schr\"{o}dinger-Virasoro Lie algebras (including both original and twisted sectors) and their deformations were respectively introduced in \cite{H1} (further investigated in \cite{H2,HU}) and \cite{RU}, in the context of non-equilibrium statistical physics, which are closely related to the Schr\"{o}dinger Lie algebra and the Virasoro Lie algebra. The Lie algebras $\LL\,(\l,\mu\in\C)$
considered in this paper called {\it deformative
Schr\"{o}dinger-Virasoro Lie algebras} (\,see \cite{RU}), possess
the  $\C$-basis $\{L_n,$ $M_n,\,Y_{s+n}\,|\,n\in\Z,\,s=0\ {\rm or}\ \frac{1}{2}\}$ with the
following non-vanishing Lie brackets:
\begin{eqnarray}\begin{array}{lll}\label{LB1}
&[\,L_n,L_m]\!\!\!&=(m-n)L_{m+n},\ \ \ \,\,
[L_n,Y_{s+m}\,\,]=(s+m-\frac{(\l+1)n}{2}+\mu)Y_{s+m+n},
\\[6pt]
&[L_n,M_m]\!\!&=(m-\l n+2\,\mu)M_{m+n},\ \ \,\, [\,Y_{s+n},Y_{s+m}\,]=(m-n)M_{2s+m+n}.\end{array}
\end{eqnarray}
There appeared a series of papers on these infinite-dimensional Lie algebras due to their interests and importance. The non-trivial vertex algebra representations of the original sector were explicitly constructed in \cite{U}. The Lie bialgebra structures, modules of intermediate series and Whittaker modules over the original Schr\"{o}dinger-Virasoro Lie algebra were respectively determined in \cite{HLS,LS,ZTL}. The derivations, central extensions and automorphism group of the extended sector were investigated in \cite{GJP2}. The automorphisms together with Verma modules and central extensions together with derivations for generalized Schr\"{o}dinger-Virasoro algebras were  respectively investigated in \cite{TZ} and \cite{WLXin}. The second cohomology group of the deformative Schr\"{o}dinger-Virasoro algebras (both original and twisted) were determined by \cite{LJ,LSZ,RU} completely. Derivations and automorphisms of twisted
deformative Schr\"{o}dinger-Virasoro Lie algebras were described in \cite{WLXu}.

For any $\C$-Lie algebra $\gg$, a symmetric bilinear forms $\varphi$ on $\gg$
is called $\gg$-invariant if it satisfies
\begin{eqnarray*}
\vp([x,y],z)=\vp(x,[y,z]),\ \ \forall\,x,y,z\in\gg.
\end{eqnarray*}
Denote $\mathrm{Inv}(\gg,\C)$ the set of all $\C$-valued symmetric $\gg$-invariant bilinear forms on the Lie algebra $\gg$. The non-degeneracy of a form $\varphi\in\mathrm{Inv}(\gg,\C)$ means $\varphi(x,\gg)=0$ implying $x=0$. It is well known that non-degenerate symmetric invariant bilinear forms paly important roles in the structures and representations of Lie algebras and the Lie algebras with nondegenerate, symmetric and invariant bilinear forms (usually called quadratic Lie algebras) play a privileged role in physics (see \cite{FS,GJP2,HPL,ZM} and their references).

The main purpose of this paper is to determine the non-degenerate invariant bilinear forms of deformative Schr\"{o}dinger-Virasoro algebras $\LL\,(\l,\mu\in\C)$
defined above, which can be formulated as the following theorem.
\begin{theo}\label{theo}
\begin{eqnarray*}
\begin{array}{lll}
&&\Inv\cong\left\{\begin{array}{lllll}
\C\varphi_{\frac{1}{2}\Z,-1}&\mbox{\rm{if}\ \ }\mu\in\frac{1}{2}\Z,\,\lambda=-1,\\[3pt]
\C\varphi_{s+\Z,-3}&\mbox{\rm{if}\ \ }\mu\in s+\Z,\,\lambda=-3,\\[3pt]
\C\varphi_{s+\Z,-5}&\mbox{\rm{if}\ \ }\mu\in s+\Z,\,\lambda=-5,\\[3pt]
\C\varphi_{s+\Z,-2}&\mbox{\rm{if}\ \ }\mu\in s+\Z,\,\lambda=-2,\\[3pt]
\C\varphi_{\frac{1}{2}\Z\DU s+\Z,-2}&\mbox{\rm{if}\ \ }\mu\in\frac{1}{2}\Z\DU s+\Z,\,\lambda=-2,\\[3pt]
0&\mbox{{\rm otherwise}},
\end{array}\right.
\end{array}
\end{eqnarray*}
where the corresponding non-vanishing components can be given as follows:
\begin{eqnarray*}
\varphi_{\frac{1}{2}\Z,-1}(L_{0},M_{-2\mu})\!\!\!&=&\!\!\! 1,\ \ \ \varphi_{\frac{1}{2}\Z,-1}(Y_{p},Y_{-p-2\mu})=-2(1-\delta_{p,-\mu}),\\[3pt]
\varphi_{s+\Z,-3}(L_{0},M_{-2\mu})\!\!\!&=&\!\!\! 1,\ \ \
\varphi_{s+\Z,-5}(L_{n},Y_{-n-\mu})=1,\\[3pt]
\varphi_{\frac{1}{2}\Z\DU s+\Z,-2}(L_{n},M_{-n-2\mu})\!\!\!&=&\!\!\! 1,\ \ \
\varphi_{s+\Z,-2}(L_{n},M_{-n-2\mu})=1,\ \ \ \varphi_{s+\Z,-2}(Y_{p},Y_{-p-2\mu})=-2,
\end{eqnarray*}
for $s=0\ {\rm or}\ \frac{1}{2}$ and all $n\in\Z$, $p\in s+\Z$.
\end{theo}

In order to present the corollary of the above theorem, we shall introduce some related notions based on a Lie algebra $\mathcal{S}$. A {\it 2-cocycle} on $\mathcal{S}$ is a $\C$-bilinear function
$\psi:\mathcal{S}\times\mathcal{S}\rar\C$ satisfying the skew-symmetry and the Jacobi identity
\begin{eqnarray}
&&\psi(x,y)=-\psi(y,x),\label{2-cocy-comm}\\[3pt]
&&\psi([x,y],z)=\psi(x,[y,z])-\psi(y,[x,z]),\label{2-cocy-Jco}
\end{eqnarray}
for any $x,y,z\in\mathcal{S}$. Denote by ${\mathcal {C}^2}(\mathcal{S},\C)$ the
vector space of 2-cocycles on $\mathcal{S}$. For any $\C$-linear
function $f:\mathcal{S}\rar\C$, one can define a 2-cocycle $\psi_f$
as follows:
\begin{eqnarray}\label{cobcoy}
\psi_f(x,y)=f([x,y]),\ \ \forall\,\,x,y\in\mathcal{S},
\end{eqnarray}
which is called a {\it 2-coboundary} or a {\it trivial 2-cocycle} on $\mathcal{S}$. Denote by ${\BB^2}(\mathcal{S},\C)$ the vector
space of 2-coboundaries on $\mathcal{S}$. A 2-cocycle $\phi$
is said to be {\it equivalent to} another 2-cocycle $\psi$
if $\phi-\psi$ is trivial. The quotient space
\begin{eqnarray}\label{dseco}
\HH^2(\mathcal{S},\C)\!=\!\mathcal {C}^2(\mathcal{S},\C)/\BB^2(\mathcal{S},\C),
\end{eqnarray}
is called the {\it second cohomology group} of $\mathcal{S}$. It is known that the second cohomology group $\HH^2(\mathcal {L}^s_{\l,\mu},\C)$ of the deformative Schr\"{o}dinger-Virasoro algebras (both original and twisted) have been determined by \cite{LJ,LSZ,RU} completely.

Any Lie algebra $\NN$ will become a {\it Leibniz algebra} if we remove the condition $[x,y]+[y,x]=0\,(\,\forall\,x,y\in\NN\,)$. Any  2-cocycle  $\psi:\NN\times\NN\rar\C$ will become a {\it Leibniz 2-cocycle} if we remove the condition \eqref{2-cocy-comm}. We use ${\mathcal {C}_L^2}(\NN,\C)$ and ${\BB_L^2}(\NN,\C)$ respectively to denote the vector spaces of Leibniz 2-cocycles and Leibniz 2-coboundaries on $\NN$. Then the corresponding quotient space $\HH_L^2(\NN,\C)\!=\!\mathcal {C}_L^2(\NN,\C)/\BB_L^2(\NN,\C)$ is called the {\it second Leibniz cohomology group} of $\NN$.

There is a natural question whether the second cohomology group is consistent with the second Leibniz cohomology group for a given Lie algebra. The non-degenerateness of all the invariant bilinear forms on the corresponding Lie algebra plays the central role. The following proposition describes the details, which can be found in some references (e.g., \cite{GJP2,HPL,P}).
\begin{prop}\label{prop}
For any Lie algebra $\frak{g}$ over $\C$, the following exact sequence holds:
\begin{eqnarray*}
0\longrightarrow\HH^2(\frak{g},\C)\stackrel{\pi}{\longrightarrow}\HH_L^2(\frak{g},\C)
\stackrel{\xi}{\longrightarrow}\mathrm{Inv}(\frak{g},\C)
\stackrel{\psi}{\longrightarrow}\HH^3(\frak{g},\C)
\end{eqnarray*}
where $\mathrm{Inv}(\frak{g},\C)$ denotes the space spanned by the symmetric invariant bilinear forms on $\frak{g}$, $\pi$ is the natural embedding map, $\xi$ is defined by $\xi(\mu)(x,y)=\mu(x,y)+\mu(y,x)$ for any $\mu\in\HH_L^2(\frak{g},\C)$, $x,y\in\frak{g}$, and $\psi$ is the Cartan-Koszul map defined by $\psi(\nu)(x,y,z)=\nu([x,y],z)$ for any $\nu\in\mathrm{Inv}(\frak{g},\C)$, $x,y,z\in\frak{g}$.
\end{prop}

The following corollary can be easily verified from Theorem \ref{theo} and Proposition \ref{prop}.
\begin{coro}\label{coro}
\begin{eqnarray*}
\begin{array}{lll}
&&\HH_L^2(\mathcal {L}^s_{\l,\mu},\C)=\left\{\begin{array}{lllll}
\HH^2(\mathcal {L}^s_{-1,\mu},\C)\oplus\C\chi_{\frac{1}{2}\Z,-1}&\mbox{\rm{if}\ \ }\mu\in\frac{1}{2}\Z,\,\lambda=-1,\\[3pt]
\HH^2(\mathcal {L}^s_{-3,\mu},\C)\oplus\C\chi_{s+\Z,-3}&\mbox{\rm{if}\ \ }\mu\in s+\Z,\,\lambda=-3,\\[3pt]
\HH^2(\mathcal {L}^s_{-5,\mu},\C)\oplus\C\chi_{s+\Z,-5}&\mbox{\rm{if}\ \ }\mu\in s+\Z,\,\lambda=-5,\\[3pt]
\HH^2(\mathcal {L}^s_{-2,\mu},\C)\oplus\C\chi_{s+\Z,-2}&\mbox{\rm{if}\ \ }\mu\in s+\Z,\,\lambda=-2,\\[3pt]
\HH^2(\mathcal {L}^s_{-2,\mu},\C)\oplus\C\chi_{\frac{1}{2}\Z\DU s+\Z,-2}&\mbox{\rm{if}\ \ }\mu\in\frac{1}{2}\Z\DU s+\Z,\,\lambda=-2,\\[3pt]
\HH^2(\mathcal {L}^s_{\l,\mu},\C)&\mbox{{\rm otherwise}},
\end{array}\right.
\end{array}
\end{eqnarray*}
where the corresponding non-vanishing components can be given as follows:
\begin{eqnarray*}
\chi_{\frac{1}{2}\Z,-1}(L_{0},M_{-2\mu})\!\!\!&=&\!\!\! 1,\ \ \ \chi_{\frac{1}{2}\Z,-1}(Y_{p},Y_{-p-2\mu})=-2(1-\delta_{p,-\mu}),\\[3pt]
\chi_{s+\Z,-3}(L_{0},M_{-2\mu})\!\!\!&=&\!\!\! 1,\ \ \
\chi_{s+\Z,-5}(L_{n},Y_{-n-\mu})=1,\\[3pt]
\chi_{\frac{1}{2}\Z\DU s+\Z,-2}(L_{n},M_{-n-2\mu})\!\!\!&=&\!\!\! 1,\ \ \
\chi_{s+\Z,-2}(L_{n},M_{-n-2\mu})=1,\ \ \ \chi_{s+\Z,-2}(Y_{p},Y_{-p-2\mu})=-2,
\end{eqnarray*}
for $s=0\ {\rm or}\ \frac{1}{2}$ and all $n\in\Z$, $p\in s+\Z$.
\end{coro}

\vs{18pt}

\cl{\bf\S2. \ Proof of Theorem \ref{theo}}\setcounter{section}{2}\setcounter{theo}{0}\setcounter{equation}{0}

\vs{12pt}

The following lemma is known from the references and also can be easily deduced.
\begin{lemm}\label{LL}
$\vp(L_m,L_n)=0,\ \,\forall\ \vp\in\Inv,\ \,m,\,n\in\Z.$
\end{lemm}

For any $m\in\Z$, $p\in s+\Z$ with $p+\mu\neq 0$, the following identities hold:
\begin{eqnarray}\label{LY01}
\vp(L_m,Y_p)=\frac{1}{p+\mu}\vp(L_m,[L_0,Y_p])
=\frac{1}{p+\mu}\vp([L_m,L_0],Y_p)
=-\frac{m}{p+\mu}\vp(L_m,Y_p),
\end{eqnarray}
which implies the following two lemmas.
\begin{lemm}\label{LY--01}
If $\mu\notin s+\Z$, then $\vp(L_m,Y_p)=0,\ \ \forall\ \vp\in\Inv,\ \,m\in\Z,\ p\in s+\Z.$
\end{lemm}

\begin{lemm}\label{LY--02}
If $\mu\in s+\Z$, then $\vp(L_m,Y_p)=0,\ \ \forall\ \vp\in\Inv,\ \,m\in\Z,\ p\in s+\Z$ with the additional conditions $p\neq-\mu$ and $m+p\neq-\mu$.
\end{lemm}

\begin{lemm}\label{LY--03}
If $\mu\in s+\Z$, then $\vp(L_m,Y_{-\mu})=0,\ \ \forall\ \vp\in\Inv,\ \,m\in\Z^*$.
\end{lemm}
{\it Proof}\ \ \
If $\mu\in s+\Z$ and $\lambda\neq-3$, $\vp(L_m,Y_{-\mu})$
can be rewritten as
\begin{eqnarray}\label{LY03}
\frac{-2}{\lambda+3}\vp(L_m,[L_1,Y_{-\mu-1}])
=\frac{-2}{\lambda+3}\vp([L_m,L_1],Y_{-\mu-1})
=\frac{2(m-1)}{\lambda+3}\vp(L_{m+1},Y_{-\mu-1}),
\end{eqnarray}
which gives
\begin{eqnarray*}
\vp(L_m,Y_{-\mu})=0,\ \ \forall\ \vp\in\Inv,\ \,m\in\Z^*,\ \mu\in s+\Z,\ \lambda\neq-3.
\end{eqnarray*}

If $\mu\in s+\Z$, $\lambda=-3$ and $m\neq 2$, $\vp(L_m,Y_{-\mu})$
can be rewritten as
\begin{eqnarray}\label{LY04}
\frac{1}{2-m}\vp([L_{m-1},L_1],Y_{-\mu})
=\frac{1}{2-m}\vp(L_{m-1},[L_1,Y_{-\mu}])
=\frac{1}{2-m}\vp(L_{m-1},Y_{1-\mu}),
\end{eqnarray}
which gives
\begin{eqnarray*}
\vp(L_m,Y_{-\mu})=0,\ \ \forall\ \vp\in\Inv,\ \,m\in\Z^*\sm\{2\},\ \mu\in s+\Z,\ \lambda=-3.
\end{eqnarray*}

If $\mu\in s+\Z$ and $\lambda=-3$, $\vp(L_2,Y_{-\mu})$ can be rewritten as
\begin{eqnarray}\label{LY05}
\frac{1}{4}\vp([L_{-1},L_3],Y_{-\mu})
=\frac{1}{4}\vp(L_{-1},[L_3,Y_{-\mu}])
=\frac{3}{4}\vp(L_{-1},Y_{3-\mu}),
\end{eqnarray}
which gives
\begin{eqnarray*}
\vp(L_2,Y_{-\mu})=0,\ \ \forall\ \vp\in\Inv,\ \mu\in s+\Z,\ \lambda=-3.
\end{eqnarray*}
Then the lemma follows.\QED

When $\mu\in s+\Z$ and $m\in\Z^*$, one has the following identities:
\begin{eqnarray}\label{LY02}
\vp(L_m,Y_{-m-\mu})=\frac{1}{m}\vp([L_0,L_m],Y_{-m-\mu})
=\frac{1}{m}\vp(L_0,[L_m,Y_{-m-\mu}])
=-\frac{\lambda+3}{2}\vp(L_0,Y_{-\mu}),
\end{eqnarray}
which combining with Lemmas \ref{LY--02} and \ref{LY--03}, concludes the following lemma.
\begin{lemm}\label{LY--04}
If $\mu\in s+\Z$, then $\vp(L_m,Y_p)=-\frac{\lambda+3}{2}\delta_{p,-m-\mu}\vp(L_0,Y_{-\mu})$, $\forall\ \vp\in\Inv$, $m\in\Z^*$, $p\in s+\Z$ and $\vp(L_0,Y_q)=0$, $\forall\,-\mu\neq q\in s+\Z$.
\end{lemm}
For any $\mu\in s+\Z$ and $m,\,n\in\Z$, the following identity
\begin{eqnarray*}
\vp([L_m,L_n],Y_{-\mu-m-n})=\vp(L_m,[L_n,Y_{-\mu-m-n}])
\end{eqnarray*}
gives
\begin{eqnarray*}
(m-n)\vp(L_{m+n},Y_{-\mu-m-n})=\big(m+\frac{\lambda+3}{2}n\big)
\vp(L_m,Y_{-\mu-m}]).
\end{eqnarray*}
If $\l\neq-3$ and $m,\,m+n\in\Z^*$, using Lemma \ref{LY--04}, one further has
\begin{eqnarray*}
\frac{\lambda+5}{2}\,n\,\vp(L_0,Y_{-\mu})=0,
\end{eqnarray*}
which implies
\begin{eqnarray*}
\vp(L_0,Y_{-\mu})=0,\ \ \forall\ \,\mu\in s+\Z\ \ {\rm for\ the \ cases\ }\ \l\neq-5,-3.
\end{eqnarray*}
Then one has the following lemma.
\begin{lemm}\label{LY--05}
For any $\mu\in s+\Z$ and $\vp\in\Inv$,
\begin{eqnarray*}
\begin{array}{lll}
&&\vp(L_m,Y_p)=\left\{\begin{array}{lll}
\delta_{p,-m-\mu}\vp(L_0,Y_{-\mu})&\mbox{\rm{if}\ \ }\l=-5,\\[3pt]
0&\mbox{\rm{if}\ \ }\l\neq-5,-3,\\[3pt]
0&\mbox{\rm{if}\ \ }\l=-3,\ \ (m,p)\neq(0,-\mu).
\end{array}\right.
\end{array}
\end{eqnarray*}
\end{lemm}

\begin{lemm}\label{LM--01}
If $\mu\notin \frac{1}{2}\Z$, then $\vp(L_n,M_m)=0$, $\forall\ \vp\in\Inv$, $m,n\in\Z$.
\end{lemm}
{\it Proof}\ \ \
For any $m,n\in\Z$ and $\mu\notin \frac{1}{2}\Z$, the following identities hold:
\begin{eqnarray}\label{LM01}
\vp(L_n,M_m)\!=\!\frac{1}{m+2\mu}\vp(L_n,[L_0,M_m])\!=\!\frac{1}{m+2\mu}\vp([L_n,L_0],M_m)
\!=\!\frac{-n}{m+2\mu}\vp(L_n,M_m),
\end{eqnarray}
which gives
\begin{eqnarray*}
\vp(L_n,M_m)=0,\ \ \forall\,\,m,n\in\Z.
\end{eqnarray*}
Then the lemma follows.\QED
\begin{lemm}\label{LM--02}
If $\mu\in\frac{1}{2}\Z$, then $\vp(L_n,M_m)=-(\lambda+1)\delta_{m,-n-2\mu}\vp(L_0,M_{-2\mu})$, $\forall\ \vp\in\Inv$, $n\in\Z^*$ and $\vp(L_0,M_m)=0$, $\forall\,-2\mu\neq m\in\Z$.
\end{lemm}
{\it Proof}\ \ \
For the case $\mu\in \frac{1}{2}\Z$, the identity (\ref{LM01}) also gives
\begin{eqnarray}\label{LMRS01}
\vp(L_n,M_m)=0,\ \ \forall\,\,m\neq-2\mu,\ \ m+n\neq -2\mu.
\end{eqnarray}
For $\lambda\neq-1$ and $n\in\Z^*$, $\vp(L_n,M_{-2\mu})$ can be rewritten as
\begin{eqnarray}\label{LM02}
\frac{-1}{\lambda+1}\vp(L_n,[L_1,M_{-2\mu-1}])
=\frac{-1}{\lambda+1}\vp([L_n,L_1],M_{-2\mu-1})
=\frac{n-1}{\lambda+1}\vp(L_{n+1},M_{-2\mu-1}),
\end{eqnarray}
which together with (\ref{LMRS01}), gives
\begin{eqnarray*}
\vp(L_n,M_{-2\mu})=0,\ \ \forall\,\,n\in\Z^*,\ \lambda\neq-1.
\end{eqnarray*}

For the case $\lambda=-1$, if $n\in\Z\sm\{0,\,2\}$, $\vp(L_n,M_{-2\mu})$ can be rewritten as
\begin{eqnarray}\label{LM03}
\frac{1}{2-n}\vp([L_{n-1},L_1],M_{-2\mu})
=\frac{1}{2-n}\vp(L_{n-1},[L_1,M_{-2\mu}])
=\frac{1}{2-n}\vp(L_{n-1},M_{1-2\mu}),
\end{eqnarray}
which together with (\ref{LMRS01}), gives
\begin{eqnarray*}
\vp(L_n,M_{-2\mu})=0,\ \ \forall\,\,n\in\Z\sm\{0,\,2\}.
\end{eqnarray*}
For the case $\lambda=-1$, $\vp(L_2,M_{-2\mu})=0$ can be obtained from the following deductions:
\begin{eqnarray*}
\vp(L_2,M_{-2\mu})=\frac{1}{4}\vp([L_{-1},L_3],M_{-2\mu})
=\frac{1}{4}\vp(L_{-1},[L_3,M_{-2\mu}])
=\frac{3}{4}\vp(L_{-1},M_{3-2\mu}).
\end{eqnarray*}
Thus
\begin{eqnarray*}
\vp(L_n,M_{-2\mu})=0,\ \ \ \forall\,\,n\in\Z^*.
\end{eqnarray*}
For any $n\in\Z^*$,  $\vp(L_n,M_{-n-2\mu})$ can be rewritten as
\begin{eqnarray}\label{LM04}
\frac{1}{n}\vp([L_0,L_n],M_{-n-2\mu})
=\frac{1}{n}\vp(L_0,[L_n,M_{-n-2\mu}])
=-(\lambda+1)\vp(L_0,M_{-2\mu}).
\end{eqnarray}
Then the lemma follows.\QED

For any $\mu\in\frac{1}{2}\Z$ and $m,\,n\in\Z$, the following identity
\begin{eqnarray*}
\vp([L_m,L_n],M_{-2\mu-m-n})\!\!&=&\!\!\vp(L_m,[L_n,M_{-2\mu-m-n}])
\end{eqnarray*}
gives
\begin{eqnarray*}
(m-n)\vp(L_{m+n},M_{-2\mu-m-n})\!\!&=&\!\!(m+n+\lambda n)\vp(L_m,M_{-2\mu-m}).
\end{eqnarray*}
If $\l\neq-1$ and $m,\,m+n\in\Z^*$, using Lemma \ref{LM--02}, one further has
\begin{eqnarray*}
(\lambda+2)\,n\,\vp(L_0,M_{-2\mu})=0,
\end{eqnarray*}
which implies
\begin{eqnarray*}
\vp(L_0,M_{-2\mu})=0,\ \ \forall\ \,\mu\in\frac{1}{2}\Z\ \ {\rm for\ the \ cases\ }\ \l\neq-2,-1.
\end{eqnarray*}
Then one has the following lemma.
\begin{lemm}\label{LM--03}
For any $\mu\in\frac{1}{2}\Z$ and $\vp\in\Inv$,
\begin{eqnarray*}
\begin{array}{lll}
&&\vp(L_n,M_m)=\left\{\begin{array}{lll}
\delta_{n,-m-2\mu}\vp(L_0,M_{-2\mu})&\mbox{\rm{if}\ \ }\l=-2,\\[3pt]
0&\mbox{\rm{if}\ \ }\l\neq-2,-1,\\[3pt]
0&\mbox{\rm{if}\ \ }\l=-1,\ \ (n,m)\neq(0,-2\mu).
\end{array}\right.
\end{array}
\end{eqnarray*}
\end{lemm}
\begin{lemm}\label{YM--01}
$\vp(Y_p,M_m)=0$, $\forall\ \vp\in\Inv$, $p\in s+\Z$, $m\in\Z$.
\end{lemm}
{\it Proof}\ \ \
If $\mu\notin s+\Z$, $\vp(Y_p,M_m)=0$ follows from the following identities:
\begin{eqnarray*}
\vp(Y_p,M_m)=\frac{\vp([L_0,Y_p],M_m])}{p+\mu}
=\frac{\vp(L_0,[Y_p,M_m])}{p+\mu}.
\end{eqnarray*}
If $\mu\in s+\Z$, $p\neq-\mu$, $\vp(Y_p,M_m)=0$.
If $\mu\in s+\Z$ and $\lambda\neq-3$, $\vp(Y_{-\mu},M_m)=0$ follows from the following identities:
\begin{eqnarray*}
\vp(Y_{-\mu},M_m)=-\frac{2\vp([L_1,Y_{-\mu-1}],M_m)}{\lambda+3}
=-\frac{2\vp(L_1,[Y_{-\mu-1},M_m])}{\lambda+3}.
\end{eqnarray*}
If $\mu\in s+\Z$, $\lambda=-3$, $s=0$ and $m\in\Z^*$, then $\mu\in\Z$, $\vp(Y_{-\mu},M_m)=0$ follows from the following identities:
\begin{eqnarray*}
\vp(Y_{-\mu},M_m)=\frac{1}{m}\vp(Y_{-\mu},[L_{-\mu},M_{m+\mu}])
=\frac{1}{m}\vp([Y_{-\mu},L_{-\mu}],M_{m+\mu})=\frac{-\mu}{m}\vp(Y_{-2\mu},M_{m+\mu}).
\end{eqnarray*}
If $\mu\neq0$, $\vp(Y_{-\mu},M_0)=0$ follows from the following identities:
\begin{eqnarray*}
\vp(Y_{-\mu},M_0)=\frac{1}{4\mu}\vp(Y_{-\mu},[L_{\mu},M_{-\mu}])
=\frac{1}{4\mu}\vp([Y_{-\mu},L_{\mu}],M_{-\mu})=\frac{-1}{4}\vp(Y_0,M_{-\mu}).
\end{eqnarray*}
If $\mu=0$, $\vp(Y_0,M_0)=0$ follows from the following identities:
\begin{eqnarray*}
\vp(Y_0,M_0)=\frac{1}{2}\vp(Y_0,[L_1,M_{-1}])
=\frac{1}{2}\vp([Y_0,L_1],M_{-1})=\frac{-1}{2}\vp(Y_1,M_{-1}).
\end{eqnarray*}
If $\mu\in s+\Z$, $\lambda=-3$, $s=\frac{1}{2}$ and $m\neq-2-2\mu$, $\vp(Y_{-\mu},M_m)=0$ follows from the following identities:
\begin{eqnarray*}
\vp(Y_{-\mu},M_m)=\frac{1}{m+2+2\mu}\vp(Y_{-\mu},[L_1,M_{m-1}])
=\frac{1}{m+2+2\mu}\vp([Y_{-\mu},L_1],M_{m-1}).
\end{eqnarray*}
If $\mu\in s+\Z$, $\lambda=-3$, $s=\frac{1}{2}$ and $m=-2-2\mu$, $\vp(Y_{-\mu},M_{-2-2\mu})=0$ follows from the following identities:
\begin{eqnarray*}
\vp(Y_{-\mu},M_{-2-2\mu})\!\!\!&=&\!\!\!-\frac{1}{4}\vp(Y_{-\mu},[L_{-1},M_{-1-2\mu}])\\
\!\!\!&=&\!\!\!-\frac{1}{4}\vp([Y_{-\mu},L_{-1}],M_{-1-2\mu})
=-\frac{1}{4}\vp([Y_{-\mu-1},M_{-1-2\mu}).
\end{eqnarray*}
Then the lemma follows.\QED
\begin{lemm}\label{MM--01}
$\vp(M_n,M_m)=0$, $\forall\ \vp\in\Inv$, $m,n\in\Z$.
\end{lemm}
{\it Proof}\ \ \ For any $m\in\Z$, $3n\neq 2s$, one has
\begin{eqnarray*}
\vp(M_n,M_m)=\frac{\vp([Y_{s-n},Y_{2n-s}],M_m)}{3n-2s}
=\frac{\vp(Y_{s-n},[Y_{2n-s},M_m])}{3n-2s}
=\frac{\vp(Y_{s-n},0)}{3n-2s}=0.
\end{eqnarray*}
For the case $s=0$ and any $m\in\Z$, one has
\begin{eqnarray*}
\vp(M_0,M_m)=\frac{1}{2}\vp([Y_{-1},Y_{1}],M_m)
=\frac{1}{2}\vp(Y_{-1},[Y_1,M_m])
=\frac{1}{2}\vp(Y_{-1},0)=0.
\end{eqnarray*}
Then the lemma follows.\QED
\begin{lemm}\label{YY--01}
If $\mu\notin\frac{1}{2}\Z$, then $\vp(Y_p,Y_q)=0$, $\forall\ \vp\in\Inv$, $p,q\in s+\Z$.
\end{lemm}
{\it Proof}\ \ \
 $\vp(Y_{s+n},Y_{s+m})$ can be written as
\begin{eqnarray}\label{YY01}
-\frac{\vp([Y_{s+n},L_0],Y_{s+m})}{s+n+\mu}
=-\frac{\vp(Y_{s+n},[L_0,Y_{s+m}])}{s+n+\mu}
=-\frac{s+m+\mu}{s+n+\mu}\vp(Y_{s+n},Y_{s+m}),
\end{eqnarray}
which implies the lemma.\QED
The identity (\ref{YY01}) also implies the following lemma.
\begin{lemm}\label{YY--02}
If $\mu\in\frac{1}{2}\Z$, then $\vp(Y_p,Y_q)=0$, $\forall\ \vp\in\Inv$, $p,q\in s+\Z$ with $p\neq-\mu$ and $p+q\neq-2\mu$.
\end{lemm}
\begin{lemm}\label{YY--03}
If $\mu\in s+\Z$, then $\vp(Y_{-\mu},Y_q)=0$, $\forall\ \vp\in\Inv$, $p,q\in s+\Z$ with $q\neq-\mu$.
\end{lemm}
{\it Proof}\ \ \
If $p=-\mu$, $\lambda\neq-3$ and $q\neq-\mu$, then
\begin{eqnarray*}
\vp(Y_{-\mu},Y_q)\!\!\!&=&\!\!\!\frac{2}{\lambda+3}\vp([Y_{-\mu-1},L_1],Y_q)
=\frac{2}{\lambda+3}\vp(Y_{-\mu-1},[L_1,Y_q])\\
\!\!\!&=&\!\!\!\frac{2(q-\frac{\lambda+1}{2}+\mu)}{\lambda+3}\vp(Y_{-\mu-1},Y_{q+1})=0.
\end{eqnarray*}
If $p=-\mu$, $\lambda=-3$ and $q\neq-\mu$, then
\begin{eqnarray*}
\vp(Y_{-\mu},Y_q)=\frac{1}{q+\mu}\vp(Y_{-\mu},[L_0,Y_q])
=\frac{1}{q+\mu}\vp([Y_{-\mu},L_0],Y_q)=0,
\end{eqnarray*}
which implies the lemma.\QED
Then Lemmas \ref{YY--02} and \ref{YY--03} imply the following lemma.
\begin{lemm}\label{YY--04}
If $\mu\in\frac{1}{2}\Z$, then $\vp(Y_p,Y_q)=0$, $\forall\ \vp\in\Inv$, $p,q\in s+\Z$ with $p+q\neq-2\mu$.
\end{lemm}
Furthermore,  for any $\vp\in\Inv$ and $\mu\in\frac{1}{2}\Z$, one has the following lemma.
\begin{lemm}\label{YY--05}
(1)\ \ If $p\neq-\mu$, then  $\vp(Y_p,Y_{-2\mu-p})=-2\vp(L_0,M_{-2\mu})$.\\
(2)\ \ If $\lambda\neq-3$, then  $\vp(Y_{-\mu},Y_{-\mu})=\frac{2(\lambda+1)}{\lambda+3}\vp(L_0,M_{-2\mu})$.
\end{lemm}
{\it Proof}\ \ \ If $\mu\in\frac{1}{2}\Z$, $p\neq-\mu$ and $p+q=-2\mu$, then
\begin{eqnarray*}
\vp(Y_p,Y_q)\!\!\!&=&\!\!\!\frac{1}{p+\mu}\vp([L_0,Y_p],Y_q)
=\frac{1}{p+\mu}\vp(L_0,[Y_p,Y_q)]\\
\!\!\!&=&\!\!\!\frac{q-p}{p+\mu}\vp(L_0,M_{-2\mu})
=-2\vp(L_0,M_{-2\mu}).
\end{eqnarray*}

If $\mu\in s+\Z$, $p=q=-\mu$, $\lambda\neq-3$, then
\begin{eqnarray*}
\vp(Y_{-\mu},Y_{-\mu})\!\!\!&=&\!\!\!\frac{-2}{\lambda+3}\vp([L_1,Y_{-\mu-1}],Y_{-\mu})
=\frac{-2}{\lambda+3}\vp(L_1,[Y_{-\mu-1},Y_{-\mu}])\\
\!\!\!&=&\!\!\!
\frac{-2}{\lambda+3}\vp(L_1,M_{-2\mu-1})
=\frac{2(\lambda+1)}{\lambda+3}\vp(L_0,M_{-2\mu}).
\end{eqnarray*}
The lemma follows.\QED

Using Lemma \ref{LM--03}, we obtain the following lemma.
\begin{lemm}\label{YY--06}
For any $\mu\in\frac{1}{2}\Z$ and $\vp\in\Inv$,
\begin{eqnarray*}
\begin{array}{lll}
&&\vp(Y_p,Y_q)=\left\{\begin{array}{lll}
-2\delta_{p,-2\mu-q}\vp(L_0,M_{-2\mu})&\mbox{\rm{if}\ \ }\l=-2,\\[3pt]
-2(1-\delta_{p,-\mu})\delta_{p,-2\mu-q}\vp(L_0,M_{-2\mu})&\mbox{\rm{if}\ \ }\l=-1,\\[3pt]
0&\mbox{\rm{if}\ \ }\l\neq-3,-2,-1,\\[3pt]
0&\mbox{\rm{if}\ \ }\l=-3,\ \ (p,q)\neq(-\mu,-\mu).
\end{array}\right.
\end{array}
\end{eqnarray*}
\end{lemm}

\ni{\it Proof of Theorem \ref{theo}}\ \ \ It can be easily verified from Lemmas
\ref{LL}, \ref{LY--01}, \ref{LY--05}, \ref{LM--01}, \ref{LM--03}---\ref{YY--01} and \ref{YY--06}.\QED

\end{document}